\documentclass[a4paper,11pt]{article}
\usepackage{amsmath,amsthm,amsfonts,amssymb,srcltx}
\usepackage{graphicx} 
\usepackage{enumerate}
\newtheorem{thm}{Theorem}

\newtheorem{rem}{Remark}

\begin{document}

\title{Hultman numbers, polygon gluings and matrix integrals}
\author{Nikita Alexeev and Peter Zograf}

\maketitle
\begin{abstract}
The Hultman numbers enumerate permutations whose cycle graph has a given number of alternating cycles (they are relevant to the Bafna-Pevzner  approach to genome comparison and genome rearrangements). We give two new interpretations of the Hultman numbers: in terms of polygon gluings and as integrals over the space of complex matrices, and derive some properties of their generating functions.
\end{abstract}
\begin{figure}[hbt]%
 \begin{center}
 \includegraphics[width=4cm]{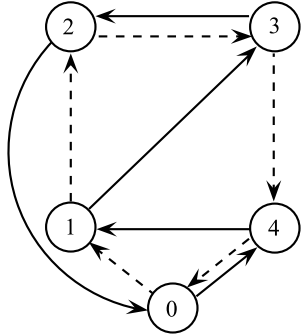}
 \caption{{\small The cycle graph $G(\pi)$ of the permutation $\pi=\genfrac{(}{)}{0pt}{}{1234}{2314},$ where the 
 grey edges are drawn by dashed arrows and the black edges are drawn by solid arrows. There are 3 alternating cycles:
 0-1-3-4-1-2-0, 2-3-2 and 4-0-4.}}
 \label{cycle_graph}%
\end{center}
\end{figure}

\noindent
{\bf Introduction.}
In the paper \cite{BP} on genome comparison and genome rearrangements, Bafna and Pevzner raised the problem of decomposing a permutation into the minimal number of ``transpositions" (here a transposition is understood as an exchange of two contiguous intervals of the permutation). An important tool they introduced to deal with this problem is the 
{\em cycle graph} of a permutation. We recall that the cycle graph of a permutation 
$\pi\in S_n$, denoted by $G(\pi)$, is the directed edge-colored graph with
vertices $\{0,1,\dots,n\}$ and edges of two colors: grey edges going from $i-1$ to $i$ and black edges 
going from $\pi_{i}$ to $\pi_{i-1},\; i=0,\ldots, n$  (throughout this note we assume that $\pi_0=0$ 
and consider $i$ modulo $n+1$). 
An {\em alternating cycle} in $G(\pi)$ is a directed cycle with edges of alternate colors.
Notice that at every vertex of $G(\pi)$ there is one incoming edge and one
outgoing edge of each color. This means that there is a unique disjoint decomposition of
the edge set of $G(\pi)$ into alternating cycles, see Fig.~\ref{cycle_graph}.

In his thesis \cite{H}, Hultman attempted to characterize the number $H(n,k)$ of permutations in $S_n$ whose cycle graph has exactly $k$ alternating cycles. These numbers, now carrying his name (see www.oeis.org/A164652), have later been studied by several authors (cf. \cite{BF} and \cite{DL} to name just few). As it is shown in \cite{BF}, the Hultman numbers are closely related to the (unsigned) Stirling numbers of the first kind $S(n,k)$ (see www.oeis.org/A008275) that count permutations in $S_{n}$ whose disjoint cycle decomposition consists of $k$ cycles:
\begin{equation}
H(n, k) = \begin{cases} 
 \frac{2\,S(n+2,k)}{(n+1)(n+2)}&\textrm{ if } n-k \textrm{ is odd}, \\ 
 \quad\quad 0 &\textrm{ otherwise.}
\end{cases}\label{stirling}
\end{equation}
A closed formula for the Hultman numbers was obtained in \cite{DL}.

In this note we give two new interpretations of the Hultman numbers in the spirit of \cite{HZ}: as numbers of certain polygon gluings and as integrals over the space of complex matrices. We also give a recursion relation for the Hultman numbers and derive some properties of their generating functions.

\begin{figure}[hbt]%
 \begin{center}
 \includegraphics[width=5cm]{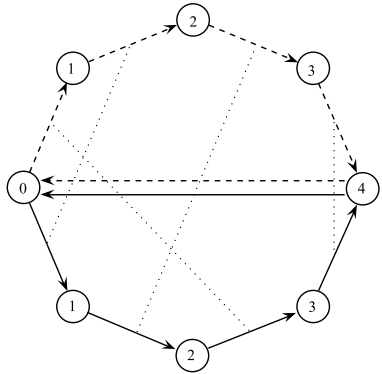}
 \caption{{\small A $2n$-gon ($n=4$) with $n$ black sides (solid arrows) and $n$ grey sides (dashed arrows). The pairs of sides that are glued together by 
 $\pi=\genfrac{(}{)}{0pt}{}{1234}{2314}$ are connected with dotted lines.}}
 \label{polygon}%
\end{center}
\end{figure}

\bigskip
\noindent
{\bf Polygon gluings.}
Consider a  $2n$-sided polygon, whose boundary consists of $n$ black sides followed by $n$ grey sides; the black sides are oriented in the counterclockwise direction, and the grey sides are oriented in the clockwise direction, see Fig.~\ref{polygon}.

Pairwise gluing of black sides with grey sides (respecting orientation) gives an orientable topological surface without boundary of topological genus $g\geq 0$ (the genus $g$ depends on the gluing). At the same time, the boundary of the polygon turns into an oriented graph with $k\geq 1$ vertices and $n$ edges. The numbers $g$ and $k$ are related by the Euler characteristic formula $2-2g=k-n+1$, so that $k=n+1-2g$. We denote by $h_g(n)$ the number of genus $g$ such gluings of a $2n$-gon.

\begin{thm} The Hultman numbers $H(n,k)$ and the numbers $h_g(n)$ of genus $g$ gluings of a $2n$-gon described above 
are related by the fomula
\begin{equation}
H(n,n+1-2g)=h_g(n).
\end{equation}
\label{bij}
\end{thm} 

\noindent 
\begin{proof}
We start with a slightly different interpretation of the cycle graph $G(\pi)$. Consider two oriented cycles (that is, 2-regular oriented graphs) of length $n+1$, one colored in grey and the other colored in black. The vertex set in both cycles is $\{0,\ldots,n\}$, but in the grey cycle the vertices follow in the clockwise order, and in the black cycle they follow in the counterclockwise order. We identify the vertex ${\pi_i}$ of the grey cycle with the vertex $i$ of the black cycle (we assume $\pi_0=0$). Obviously, the obtained graph coincides with the cycle graph $G(\pi)$, see Fig.~\ref{cycle_graph}.

We label the black sides of the polygon by numbers from 1 to $n$ in the counterclockwise order, and the grey sides by numbers from 1 to $n$ in the clockwise order, both times starting from the initial vertex $0$. Clearly, a gluing of a $2n$-gon of the type considered above is uniquely described by a permutation $\pi \in S_n$, where $\pi_i$ is the number of the grey side identified with the $i$th black side. Let us cut the polygon along the diagonal $(n,0)$, i.e.,
we add one black edge and one grey edge connecting the vertex $n$ to the vertex 0, see Fig.~\ref{polygon}.
Now we have two $n$-gons, one with black boundary and the other with grey boundary, whose sides are pairwise identified by means of the permutation $\pi$ ($\pi_0=0$). These two boundaries glued together give a graph that we denote by $\Gamma(\pi)$. The construction is quite similar to that of the cycle graph $G(\pi)$, but instead of gluing vetices we now glue edges
according to the same rule. 
The graphs $G(\pi)$ and $\Gamma(\pi)$ are closely related to each other: it is straightforward to verify that there is a one-to-one correspondence between the alternating cycles in the cycle graph $G(\pi)$ and the vertices in the polygon gluing graph $\Gamma(\pi)$.
To complete the proof, we recall that $k=n+1-2g$, where $k$ is the number of vertices of $\Gamma(\pi)$, and $g$ is the genus of the glued surface.
\end{proof}

\bigskip
\noindent
{\bf Matrix integral.}
Denote by $M(N)={\rm Mat}_{\mathbb{C}}(N\times N)$ the linear space of complex $N\times N$ matrices; the (complex) dimension of $M(N)$ is $N^2$. The space $M(N)$ has a natural Gaussian probabilistic measure
\begin{equation}
d\mu_N=\left(\frac{1}{2\pi \sqrt{-1}}\right)^{N^2} e^{-{\rm Tr} (XX^*)} \bigwedge_{i,j=1}^N dx_{ij}\wedge d\bar{x}_{ij},
\end{equation}
where $X=\{x_{ij}\}_{i,j=1}^N \in M(N)$, the star  $*$ denotes the Hermitian conjugation and Tr is the trace. Note that the space $M(N)$ equipped with the measure $\mu_N$ is also called the complex Ginibre ensemble. 

\begin{thm}
Put
\begin{equation}
p_n(N)=\sum_{g=0}^{[n/2]} H(n,n+1-2g)\,N^{n-2g+1},\label{polynom}
\end{equation}
where $H(n,k)$ are the Hultman numbers. Then
\begin{equation}
p_n(N)=\int_{M(N)}{\rm Tr}(X^nX^{*n})\,d\mu_N.\label{matrix}
\end{equation} 

\label{matrix_integral}
\end{thm}
\begin{rem}
{\rm More general matrix integrals over the space $M(N)$ are considered in \cite{AGT}.}
\end{rem}
\begin{rem}
{\rm Below is a list of the several first polynomials $p_n(N)$:}
\begin{align*}
&p_0(N)=N,\\
&p_1(N)=N^2,\\
&p_2(N)=N^3 + N,\\
&p_3(N)=N^4 + 5 N^2,\\
&p_4(N)=N^5 + 15 N^3 + 8 N,\\
&p_5(N)=N^6 + 35 N^4 + 84 N^2,\\
&p_6(N)=N^7 + 70 N^5 + 469 N^3 + 180 N,\\
&p_7(N)=N^8  + 126 N^6  + 1869 N^4  + 3044 N^2,\\
&p_8(N)=N^9  + 210 N^7  + 5985 N^5  + 26060 N^3  + 8064 N,\\
&p_9(N)=N^{10} + 330 N^8 + 16401 N^6  + 152900 N^4  + 193248 N^2.  
\end{align*}
\end{rem}
\noindent
\begin{proof}
It is a fairly standard exercise in t'Hooft graphic calculus to reduce the matrix integral in Eq.~(\ref{matrix}) to a sum over Feynman diagrams (polygon gluings), cf. e.g. \cite{M}, \cite{Z}. We will briefly explain how it works. By definition we have
$$
{\rm Tr}(X^nX^{*n})=\sum_{i_1=1}^N\ldots\sum_{i_{2n}=1}^N x_{i_1i_2}\ldots x_{i_ni_{n+1}}\bar{x}_{i_{1}i_{2n}}\ldots\bar{x}_{i_{n+2}i_{n+1}},
$$
and a simple computation shows that
$$
\int_{M(N)}x_{ij}\bar{x}_{kl}d\mu_N=\delta_{ik}\delta_{jl},\qquad\int_{M(N)}x_{ij}x_{kl}d\mu_N=\int_{M(N)}\bar{x}_{ij}\bar{x}_{kl}d\mu_N=0.
$$
Applying Wick's formula (cf. \cite{M}, \cite{Z}), we get
\begin{align*}
&\int_{M(N)}x_{i_1i_2}\ldots x_{i_ni_{n+1}}\bar{x}_{i_{1}i_{2n}}\ldots\bar{x}_{i_{n+2}i_{n+1}}d\mu_N\\
&=\sum_{\pi\in S_n}\int_{M(N)}x_{i_1i_2}\bar{x}_{i_{\alpha_1+1}i_{\alpha_1}}d\mu_N\cdots\int_{M(N)}x_{i_ni_{n+1}}\bar{x}_{i_{\alpha_n+1}i_{\alpha_n}}d\mu_N\\
&=\sum_{\pi\in S_n}\delta_{i_1 i_{\alpha_1+1}}\delta_{i_2 i_{\alpha_1}}\cdots\delta_{i_ni_{\alpha_n+1}}\delta_{i_{n+1}i_{\alpha_n}},
\end{align*}
where $\alpha_j=2n+1-\pi_j$ {(we assume that $i_{2n+1}=i_1$)}. Therefore,
\begin{align*}
\int_{M(N)}{\rm Tr}(X^nX^{*n})\,d\mu_N=\sum_{\pi\in S_n}\sum_{i_1=1}^N\ldots\sum_{i_{2n}=1}^N
\delta_{i_1 i_{\alpha_1+1}}\delta_{i_2 i_{\alpha_1}}\cdots\delta_{i_ni_{\alpha_n+1}}\delta_{i_{n+1}i_{\alpha_n}}.
\end{align*}
We note that the pairs of indices $\{i_k i_{k+1}\}$ correspond to the black edges of the polygon on Fig.~(\ref{polygon}), and the pairs of indices 
$\{i_{\alpha_k+1} i_{\alpha_k}\}$ correspond to the grey edges, so there is a one-to one correspondence between the pairings of indices and polygon gluings. Moreover, it is not hard to see that for a given $\pi\in S_N$
$$
\sum_{i_1=1}^N\ldots\sum_{i_{2n}=1}^N
\delta_{i_1 i_{\alpha_1+1}}\delta_{i_2 i_{\alpha_1}}\cdots\delta_{i_ni_{\alpha_n+1}}\delta_{i_{n+1}i_{\alpha_n}}=N^{n-2g+1},
$$
where $g$ denotes the genus of the surface glued from the $2n$-gon by means of $\pi$. This yields
$$
\int_{M(N)}{\rm Tr}(X^nX^{*n})\,d\mu_N=\sum_{g=0}^{[n/2]}h_g(n)N^{n-2g+1},
$$
and Eq.~(\ref{matrix}) now follows from Theorem \ref{bij}.\end{proof}

\bigskip
\noindent
{\bf Generating functions and recursions.}
Here we collect some simple facts about the recursive relations and generating functions for the Hultman numbers that we did not find in the literature. 

Consider the generating functions
\begin{equation}
F(x,N)=\sum_{g=0}^\infty\sum_{n=2g}^\infty H(n,n+1-2g)N^{n-2g+1}\frac{x^n}{n!}
\end{equation}
and
\begin{equation}
H_g(x)=\sum_{n=2g}^\infty H(n,n+1-2g)x^n.
\end{equation}

\begin{thm}
We have
\begin{enumerate}[{(i)}]
\item{$$F(x,N)=\frac{1}{x^2}\left(\frac{1}{(1-x)^N}-(1+x)^N\right);$$}
\item{$H(n,n+1-2g)=h_g(n)$ satisfy the recursion
$$(n+2)h_g(n)=(2n+1)h_g(n-1)-(n-1)h_g(n-2)+n^2(n-1)h_{g-1}(n-2);$$}
\item{the polynomials $p_n(N)$ defined by Eq.~(\ref{polynom}) satisfy the recursion
$$(n+2)p_n(N)=(2n+1)Np_{n-1}(N)+(n-1)(n^2-N^2)p_{n-2}(N)$$
with $p_0=N,\,p_1=N^2;$}
\item{$$H_0(x)=\frac{1}{1-x}\,,\quad\quad H_g(x)=\frac{P_g(x)}{(1-x)^{1+4g}}\,,\quad g\geq 1,$$
where $P_g(x)=\sum_{i=2g}^{4g-2}a_{g,i}x^i$ is a polynomial with integer coefficients,
$a_{g,2g}=\frac{(2g)!}{g+1},\, a_{g,4g-2}=1$, and  $P_g(1)=\frac{(4g-1)!!}{2g+1}$.}
\end{enumerate}
\end{thm}
\begin{rem}
{\rm Several first polynomials $P_g(x)$ are listed below:}
\begin{align*}
P_0(x)=&1,\\
P_1(x)=&x^2,\\
P_2(x)=&x^4(8+12x+x^2),\\
P_3(x)=&x^6(180+704x+528x^2+72x^3+x^4),\\
P_4(x)=&x^8(8064+56160x+98124x^2+53792x^3+8760x^4+324x^5+x^6),\\
P_5(x)=&x^{10}(604800+6356160x+19083456x^2+21676144x^3+\\
&\qquad\qquad +9936360x^4+1759520x^5+103040x^6+1344x^7+x^8).
\end{align*}
{\rm Remarkably, all polynomials $P_g(x)$ have positive integer coefficients. Moreover, the integers $P_g(1)$ are well known
(see www.http://oeis.org/A035319) -- they enumerate genus $g$ orientable gluings of a $2g$-gon \cite{HZ}, or the permutations in
$S_{4g-1}$ whose cycle graph alternating cycles are all of length 2 \cite{DL}.}
\end{rem}
\noindent
\begin{proof} Part $(i)$ follows from Eq.~(\ref{stirling}) and the fact that
$$
(1+x)^N=\sum_{n=0}^{\infty}\sum_{k=0}^n (-1)^{n+k}S(n,k) N^k \frac{x^n}{n!},
$$
where $S(n,k)$ are the Stirling numbers of the first kind. 
Similarly, the recursion $S(n+1,k)=S(n,k-1)+n\,S(n,k)$ for the Stirling numbers immediately implies $(ii)$. Part $(iii)$ is a direct consequence of 
$(ii)$. The proof of $(iv)$ is by induction on $g$ and follows the proof of Theorem 1 in \cite{APRW}. 
The cases $g=0,1$ being easy,
assume that the statements of part $(iv)$ of the theorem hold for $g-1,\; g\geq 2$. Put $\tilde{H}_g(x)=x^2 H_g(x)$,
then the recursion $(ii)$ is equivalent to the ODE
$$
(1-x)^2 \tilde{H}'_g(x) +(1-x)\tilde{H}_g(x)=x^4\tilde{H}_{g-1}'''(x)+2x^3\tilde{H}_{g-1}''(x) 
$$
with initial condition $\tilde{H}_g(0)=0$. Therefore, we have
\begin{equation}
\tilde{H}_g(x)=(1-x)\int_0^x\frac{t^4\tilde{H}_{g-1}'''(t)+2t^3\tilde{H}_{g-1}''(t)}{(1-t)^3}dt.\label{int}
\end{equation}
The elementary formula
\begin{equation}
\left(\frac{x^\alpha}{(1-x)^\beta}\right)'
=\frac{\alpha x^{\alpha-1}+(\beta-\alpha)x^\alpha}{(1-x)^{\beta+1}}\label{ef}
\end{equation}
immediately yields
\begin{align}
x^4\left(\frac{x^\alpha}{(1-x)^\beta}\right)'''&+\;2x^3\left(\frac{x^\alpha}{(1-x)^\beta}\right)''\nonumber\\
&=\frac{\alpha^2(\alpha-1)x^{\alpha+1}+\dots+(\beta-\alpha)^2(\beta-\alpha+1)x^{\alpha+4}}{(1-x)^{\beta+3}}.\label{elem}
\end{align}
Since, by assumption,
$$\tilde{H}_{g-1}(x)=\frac{x^2 P_{g-1}(x)}{(1-x)^{4g-3}}=\frac{\sum_{i=2g-2}^{4g-6}a_{g-1,i}\,x^{i+2}}{(1-x)^{4g-3}},$$
applying Eq.~(\ref{elem}) we get that
\begin{equation}
\frac{x^4\tilde{H}_{g-1}'''(x)+2x^3\tilde{H}_{g-1}''(x)}{(1-x)^{3}}=\frac{Q_g(x)}{(1-x)^{4g+3}},\label{eqq}
\end{equation}
where $Q_g(x)=\sum_{i=2g+1}^{4g}q_{g,i}x^i$ is a polynomial with integer coefficients,
\begin{align*}
&q_{g,2g+1}=(2g)^2(2g-1)a_{g-1,2g-2}=2(2g)!,\\
&q_{g,4g}=2a_{g-1,4g-6}=2.
\end{align*}
Consider the Laurent expansion
\begin{equation}
\frac{Q_g(x)}{(1-x)^{4g+3}}=\sum_{i=3}^{4g+3}\frac{r_{g,i}}{(1-x)^i},\label{laurent}
\end{equation}
then we have
$$\frac{\tilde{H}_g(x)}{1-x}=\sum_{i=2}^{4g+2}\frac{r_{g,i+1}}{i(1-x)^{i}}+C,$$
where the initial condition $\tilde{H}_g(0)=0$ implies that
$$C=-\sum_{i=2}^{4g+2}\frac{r_{g,i+1}}{i}.$$
Now put
\begin{equation}
\tilde{P}_g(z)=\sum_{i=2}^{4g+2}\frac{r_{g,i+1}}{i}((1-x)^{4g+2-i}-(1-x)^{4g+2})
=\sum_{i=0}^{4g+2}p_{g,i}x^i.\label{P_g}
\end{equation}
By construction, we have $p_{g,0}=0$, therefore $\tilde{H}_g(x)=\tilde{P}_g(x)/(1-x)^{4g+1}$ since they both satisfy the same first order ODE with the same initial condition. Moreover, since
$h_{g}(1)=\ldots=h_{g}(2g-1)=0$, we also have $p_{g,1}=\ldots=p_{g,2g+1}=0$. Inverting (\ref{ef}), we
see that 
\begin{align*}
&a_{g,2g}=p_{g,2g+2}=q_{g,2g+1}/(2g+2)=(2g)!/(g+1),\\
&a_{g,4g-2}=p_{g,4g}=q_{g,4g}/2=1
\end{align*}
as claimed. We also see that $P_g(x)=\tilde{P}_g(x)/x^2=(1-x)^{4g+1}H_g(x)$ must have integral coefficients because
$H_g(x)$ does.

To complete the proof it is sufficient to show that $$P_g(1)=\frac{(4g-1)(4g-3)(2g-1)}{2g+1}\,P_{g-1}(1)$$ 
(note that $P_0(1)=P_1(1)=1)$.
We have
\begin{align*}
&\tilde{H}'_{g-1}(x)=\frac{(1-x)\tilde{P}'_{g-1}(x)+(4g-3)\tilde{P}_{g-1}(x)}{(1-x)^{4g-2}}=\frac{P_{g,1}(x)}{(1-x)^{4g-2}},\\
&\tilde{H}''_{g-1}(x)=\frac{(1-x)P'_{g,1}(x)+(4g-2)P_{g,1}(x)}{(1-x)^{4g-1}}=\frac{P_{g,2}(x)}{(1-x)^{4g-1}},\\
&\tilde{H}'''_{g-1}(x)=\frac{(1-x)P'_{g,2}(x)+(4g-1)P_{g,2}(x)}{(1-x)^{4g}},
\end{align*}
and from Eq.~(\ref{eqq}) it then follows that
$$Q_g(x)=(1-x)(x^4P'_{g,2}(x)+2x^3P_{g,2}(x))+(4g-1)x^4P_{g,2}(x).$$
From here we easily get
\begin{align*}
&P_{g,1}(1)=(4g-3)P_{g-1}(1),\\ 
&P_{g,2}(1)=(4g-2)P_{g,1}(1)=(4g-2)(4g-3)P_{g-1}(1),\\ 
&Q_{g,1}(1)=(4g-1)P_{g,2}(1)=(4g-1)(4g-2)(4g-3)P_{g-1}(1).
\end{align*}
Clearly, $Q_{g,1}(1)=r_{g,4g+3}$ in the Laurent
expansion (\ref{laurent}), and from Eq.~(\ref{P_g}) we obtain 
$P_g(1)=\frac{1}{4g+2}\,Q_{g,1}(1)=\frac{(4g-1)(4g-2)(4g-3)}{4g+2}\,P_{g-1}(1)$ as claimed.
\end{proof}
\bigskip
\noindent
{\bf Acknowledgements.} 
The work of NA is supported in part by the Chebyshev Laboratory (Department of Mathematics and Mechanics, St. Petersburg State University) under the Russian Government grant 11.G34.31.0026, by the Federal target programme ``Research and Pedagogical Cadre for Innovative Russia" 2010-1.1-111-128-033 and by the RFBR grant 11-01-00310-a.
The work of PZ is partially supported by the RFBR grants 11-01-12092-OFI-M-2011 and 11-01-00677-a, and he
gratefully acknowledges the hospitality and support of QGM (Aarhus) and SCGP (Stony Brook).

\end{document}